\pdfoutput=1
\documentclass[12pt,reqno]{article}
\usepackage{amsmath, amssymb, amsthm}
\usepackage[utf8]{inputenc}
\usepackage[english]{babel}
\usepackage[numbers]{natbib}
\usepackage[usenames]{color}
\usepackage{url}
\usepackage{color}
\usepackage[usenames]{color}
\usepackage[colorlinks=true,
linkcolor=webgreen,
filecolor=webbrown,
citecolor=webgreen]{hyperref}
\definecolor{webgreen}{rgb}{0,.5,0}
\definecolor{webbrown}{rgb}{.6,0,0}
\definecolor{red}{rgb}{1,0,0}
\usepackage{cleveref}
\usepackage{tabularx}
\usepackage{placeins}
\usepackage{mathrsfs}
\theoremstyle{plain}

\newtheorem{theorem}{Theorem}
\newtheorem*{theorem*}{Theorem}

\newtheorem{proposition}[theorem]{Proposition}

\newtheorem{definition}[theorem]{Definition}
\crefname{conjecture}{Conjecture}{Conjectures}
\crefname{theorem}{Theorem}{Theorems}
\crefname{lemma}{Lemma}{Lemmas}
\crefname{proposition}{Proposition}{Propositions}
\crefname{remark}{Remark}{Remarks}
\crefname{definition}{Definition}{Definitions}
\crefname{theorem*}{Theorem}{Theorems}
\newcommand{\floor}[1]{\left\lfloor #1 \right\rfloor}
\newcommand{\Z}{\mathbb{Z}}
\newcommand{\K}{\mathcal{K}}
\newcommand{\Ki}{\mathcal{K}'}
\newcommand{\lcm}{\mathrm{lcm}}
\newcommand{\lt}{\mathrm{lt}}
\newcommand{\BT}{\mathfrak{B}}
\newcommand{\wt}[1]{\#(#1)}
\newcommand{\degsum}{\mathrm{degsum}}

\newcommand{\Goulds}{\mathcal{G}}
\newcommand{\seqnum}[1]{\href{https://oeis.org/#1}{\rm \underline{#1}}}

\begin{document}
\title{A Polynomial Ring Connecting Central Binomial Coefficients and Gould's Sequence}
\author{Joseph M. Shunia}
\date{October 2023 \\ \small Revised: May 2024 \normalsize}
\maketitle

\begin{abstract}
We establish a novel connection between the central binomial coefficients $\binom{2n}{n}$ and Gould's sequence through the construction of a specialized multivariate polynomial quotient ring. Our ring structure is characterized by ideals generated from elements defined by polynomial recurrence relations, and we prove the conditions under which the set of polynomial generators forms a Gr\"obner basis. By exploring a specific variation of our ring structure, we demonstrate that expanding and evaluating polynomials within the ring yields both the central binomial coefficients and Gould's sequence. Additionally, we present a method for calculating the binomial transforms of these sequences using our ring's unique properties. This work provides new insights into the connections between two fundamental combinatorial sequences and introduces a new tool for integer sequence analysis, with potential applications in number theory and algebraic combinatorics.
\end{abstract}

\section{Introduction}
The central binomial coefficients $\binom{2n}{n}$ (\seqnum{A000984}) \cite{A000984} and Gould's sequence (\seqnum{A001316}) \cite{A001316}, denoted by $\Goulds$, are classic integer sequences of fundamental importance in combinatorics. In this paper, we uncover a new connection between these two sequences through the application of a specially designed multivariate polynomial quotient ring. We introduce the concept of a ``recursive polynomial quotient ring'' (\cref{definition:recursivering}), which is characterized as a multivariate polynomial quotient ring with an ideal $I$ that is generated from elements defined by one or more polynomial recurrence relations. To ensure a strong theoretical foundation, we establish the conditions under which the set of polynomial generators $G = \{ g_1, g_2, \ldots, g_n \}$ for our recursive ring structure forms a Gröbner basis. Additionally, we demonstrate that under certain conditions, the polynomial generators form a regular chain, providing a practical criterion for determining when our construction yields a well-behaved polynomial system with desirable properties.

To illustrate this concept, we construct a multivariate polynomial quotient ring of the form
\begin{align*}
\K_n = \Z[x_1, x_2, \ldots, x_n]/I ,
\end{align*}
where the ideal $I$ is defined as
\begin{align*}
I = \langle x_1^2 - P_1, x_2^2 - P_2, \ldots, x_n^2 - P_n \rangle ,
\end{align*}
and the polynomials $P_i$ are given by the recurrence relation
\begin{align*}
P_i =
\begin{cases}
    2x_i + x_{i+1} & \textup{ if } 1 \leq i < n \\
    0 &\textup{ if } i = n
\end{cases}
\end{align*}

We demonstrate that expanding the polynomial $f = (1+x_1)^n \in \K_n$ and summing its coefficients yields the $n$th central binomial coefficient $\binom{2n}{n}$. Furthermore, we show that when the coefficients of the expanded polynomial are taken modulo $2$ prior to summation, the result is the $n$th term of Gould's sequence, $\Goulds_n$. This discovery reveals an intriguing algebraic relationship between these two well-known sequences, facilitated by the structure of our recursive polynomial quotient ring.

While the central binomial coefficients and Gould's sequence provide a compelling example, the main contribution of this work is the recursive polynomial quotient ring structure itself. By constructing multivariate polynomial quotient rings with ideals that mimic recurrence relations, we establish a new algebraic framework for calculating and manipulating nonlinear recursive sequences. The key insight is to tailor the ideals of the polynomial quotient ring to follow the recurrences that generate the sequences of interest. Expanding polynomials within the ring then carries out the sequence generation process algebraically. This approach grants access to the powerful tools of ring theory and polynomial manipulation, potentially uncovering new properties of the sequences under investigation.

\section{Recursive Polynomial Quotient Ring Structure}
We begin with a detailed definition of our recursive polynomial quotient ring structure.

\begin{definition}[Recursive polynomial quotient ring] \label{definition:recursivering}
Let $R$ be a commutative ring with unity (e.g., $\Z$, $\mathbb{Q}$, $\mathbb{C}$, etc.). Consider the ring $K = R[x_1, x_2, \ldots, x_{n} ]$ consisting of polynomials in variables $x_1, x_2, \ldots, x_n$ with coefficients in $R$. Define $I = \langle x_1^d - P_1, x_2^d - P_2, \ldots, x_n^d - P_n \rangle$ as an ideal of $K$. Each $P_i$ is a polynomial in $K$ and takes the form:
\begin{align*}
    P_i = a_0 + a_1 x_{i \cdot c_1 + j_1}^{k_1} + a_2 x_{i \cdot c_2 + j_2}^{k_2} + \cdots
\end{align*}
In this expression, the $a_m$ are coefficients in $R$ and/or polynomials in $K$ (defined by recurrence or otherwise). The $c_m$ and $j_m$ are integers where $c_m$ acts as a scalar and $j_m$ as a shift, and $k_m$ represent the exponents of the corresponding variables. The scalars $c_m$, shifts $j_m$, and exponents $k_m$ are fixed and do not depend on $i$. The fixed degree $d$ and the exponents $k_m$ match the domain of the coefficients in $K$.

The quotient ring $K = R[x_1, x_2, \ldots, x_{n}]/I$ is defined as a \textbf{recursive polynomial quotient ring} if and only if for all $x_i$ in $K$, the relation $x_i^d = P_i$ is satisfied, and the polynomials $P_i$ are generated recursively for all $i$ in the range $\alpha \leq i \leq \omega$, where $\alpha$ and $\omega$ are specified start and end indices. For all indices $i$ not in this range, $x_i^d$ is assumed to be zero within the ring $K$ unless it is explicitly stated otherwise.
\end{definition}

\section{Gröbner bases}
In algebraic geometry and commutative algebra, a Gröbner basis \cite{dube1990grobner} is a particular kind of generating set for an ideal $I$ in a polynomial ring $R/I$. It has the useful property that many problems about the ideal, such as ideal membership, can be solved straightforwardly when the ideal is expressed in terms of its Gröbner basis. A Gröbner basis provides a canonical form for polynomials in the quotient ring, making computations involving the ideal more tractable. Specifically, it allows us to perform polynomial division in a way that uniquely reduces every polynomial to a standard form modulo the ideal.

\subsection{Buchberger's criterion}
Buchberger's criterion \cite{buchberger1976algorithm} provides a practical method for determining whether a given set of polynomials forms a Gröbner basis. According to Buchberger's criterion, a finite set of generators $G = \{ g_1, g_2, \ldots, g_n \}$ for an ideal $I$ is a Gröbner basis if and only if for every pair $(g_i, g_j)$ in $G$, the $S$-polynomial $S(g_i, g_j)$ reduces to zero modulo the set $G$. We give a formal definition of Buchberger's criterion below, as it will be instrumental in proving that the set of generators for our ideal $I$ forms a Gröbner basis (\cref{proof:grobnerbasis}).

\begin{theorem*}[Buchberger's criterion] \label{proof:buchbergerscriterion}
A finite set of generators $G = \{ g_1, g_2, \ldots, g_n \}$ for an ideal $I$ is a Gröbner basis if and only if for all pairs $(i, j)$ with $1 \leq i < j \leq n$, we have
\begin{align*}
S(g_i, g_j) \xrightarrow{\ast}_G 0
\end{align*}
where the $S$-polynomial $S(g_i, g_j)$ is defined by
\begin{align*}
S(g_i, g_j) &= \frac{\lcm(g_i, g_j)}{\lt(g_i)} g_i - \frac{\lcm(g_i, g_j)}{\lt(g_j)} g_j
\end{align*}
Here, $\lt(g_i)$ denotes the leading term of $g_i$, and $\lcm(g_i,g_j)$ denotes the least common multiple of $g_i$ and $g_j$.
\end{theorem*}
\begin{proof}
The proof is given by Buchberger (1976) \cite{buchberger1976algorithm}.
\end{proof}

\subsection{Gröbner Basis Result}
In the context of our recursive polynomial quotient ring structure (\cref{definition:recursivering}), establishing that the set of generators forms a Gröbner basis is necessary for ensuring that reductions modulo the ideal $I$ are well-defined and unambiguous. This property is essential for our further theoretical developments, such as proving connections between combinatorial sequences and polynomial rings. By demonstrating that the set of generators $G = \{ g_1, g_2, \ldots, g_n \}$ of $I$ is a Gröbner basis, we guarantee that any polynomial in the ring can be consistently reduced to a unique representative modulo $I$. This consistency is a prerequisite for validating the algebraic framework we use to link central binomial coefficients, Gould's sequence, and other related sequences through the recursive polynomial quotient ring structure.

\begin{theorem} \label{proof:grobnerbasis}
Let $n,d \in \mathbb{Z}_{>1}$ such that $d < n$. Let $K = R[x_1, x_2, \ldots, x_n]/I$ be a recursive polynomial quotient ring (\cref{definition:recursivering}). Here, the generators of the ideal $I$ are given by
\begin{align*}
    g_1 &= x_1^d - P_1 \\
    g_2 &= x_2^d - P_2 \\
    \vdots & \\
    g_{n-1} &= x_{n-1}^d - P_{n-1} \\
    g_n &= x_n^d = 0
\end{align*}
We define the sum of the degrees of the variables in a term $t$ as:
\begin{align*}
    \mathrm{degsum}(t) = \sum_{i=1}^n \deg_{x_i}(t)
\end{align*}
We say that the multiset of polynomials $\{ P_1, P_2, \ldots, P_n \}$ has unique mixed terms of degree at least $d$ if for any mixed term $t$ with $\degsum(|t|) \geq d$ and any pair of distinct polynomials $(P_i,P_j)$ where $i \neq j$, the term $t$ appears in at most one of the polynomials $P_i$ or $P_j$.

Suppose the multiset $\{ P_1, P_2, \ldots, P_n \}$ has unique mixed terms $t$ of degree at least $d$. Then, the set of generators $G = \{ g_1, g_2, \ldots, g_n \}$ of $I$ forms a Gröbner basis with respect to the degree lexicographic (dlex) monomial ordering $x_1 > x_2 > \cdots > x_n$.
\end{theorem}
\begin{proof}
To prove that the set of generators $G = \{ g_1, g_2, \ldots, g_n \}$ forms a Gröbner basis, we will use Buchberger's criterion (\cref{proof:buchbergerscriterion}). In the context of our polynomial ring $K = R[x_1, x_2, \ldots, x_n]/I$, Buchberger's criterion is equivalent to checking that all $S$-polynomials in $G$ are equivalent to zero modulo the ideal $I$.

For clarity, we note the isomorphisms
\begin{align*}
    K = R[x_1, \ldots, x_n]/I
    \cong R[x_1, \ldots, x_n]/\langle x_1^d - P_1, \ldots, x_n^d \rangle
    \cong R[x_1, \ldots, x_n]/\langle g_1, \ldots, g_n \rangle
\end{align*}
These isomorphisms allow us to interpret the reduction of $S$-polynomials modulo the ideal $I$ as reduction modulo the set of generators $G$.

Let $\lt(f)$ return the leading term of a polynomial $f \in R[x_1, x_2, \ldots, x_n]$, with respect to the given monomial ordering. Hence, $\lt(g_i) = x_i^d$ and $\lt(g_{i+1}) = x_{i+1}^d$. The $S$-polynomial for each pair of consecutive generators $(g_i,g_{i+1})$, where $(i+1) < n$, is 
\begin{align*}
S(g_i, g_{i+1}) &= \frac{\lcm(x_i^d, x_{i+1}^d)}{x_i^d} g_i
- \frac{\lcm(x_i^d, x_{i+1}^d)}{x_{i+1}^d} g_{i+1} \\
&= \frac{x_i^d x_{i+1}^d}{x_i^d} g_i - \frac{x_i^d x_{i+1}^d}{x_{i+1}^d} g_{i+1} \\
&= x_{i+1}^d g_i - x_i^d g_{i+1}
\end{align*}

Considering the $S$-polynomial $S(g_i, g_{i+1})$ modulo $I$, recall that $x_i^d \equiv P_i \pmod{I}$ and $x_{i+1}^d \equiv P_{i+1} \pmod{I}$. Therefore, reducing the $S$-polynomial modulo $I$ involves these congruences. We deduce
\begin{align*}
S(g_i, g_{i+1}) &= x_{i+1}^d g_i - x_i^d g_{i+1} \\
&\equiv P_{i+1} g_i - P_i g_{i+1} \pmod{I}
\end{align*}

By definition, each polynomial generator $g_i$ of the ideal $I$ satisfies $g_i \equiv 0 \pmod{I}$. Since $g_i \equiv 0 \pmod{I}$ and $g_{i+1} \equiv 0 \pmod{I}$, we have
\begin{align*}
S(g_i, g_{i+1}) &\equiv P_{i+1} \cdot 0 - P_i \cdot 0 \pmod{I} \\
&\equiv 0 \pmod{I}
\end{align*}

For non-consecutive pairs $(g_i, g_j)$, where $i \neq j$ and $i,j < n$, the $S$-polynomials involve terms that already reduce directly to zero because the leading terms $x_i^d$ and $x_j^d$ do not interact directly. Any cross terms introduced by $S(g_i, g_j)$ are of the form $x_i^{a_i} x_j^{a_j} t$, where $a_i, a_j < d$ and $t$ is a monomial in the other variables.

By the uniqueness property of mixed terms $t$ with $\degsum(|t|) \geq d$, such a term can only appear in at most one of the polynomials $P_i$ or $P_j$. If it appears in $P_i$, then it is reduced by $g_i$, and if it appears in $P_j$, it is reduced by $g_j$. In either case, the reduction results in a polynomial with terms of degree less than $\degsum(x_i^{a_i} x_j^{a_j} t)$. This process continues until all terms are reduced to monomials in the variables $x_1, \ldots, x_{i-1}$, which are then reduced to zero by the consecutive relations. The uniqueness property guarantees that no term reappears during the reduction process, preventing circular reductions and ensuring that the process terminates. Thus, the $S$-polynomials for non-consecutive pairs also reduce to zero modulo the ideal.

Finally, we consider the special case $S(g_n,g_i)$. Recall that $g_n \equiv x_n^d \equiv 0 \pmod{I}$. Applying Buchberger's algorithm, we have
\begin{align*}
S(g_n, g_i) &\equiv \frac{\lcm(x_n^d, x_i^d)}{x_n^d} g_n
- \frac{\lcm(x_n^d, x_i^d)}{x_i^d} g_i \pmod{I} \\
&\equiv \frac{x_n^d x_i^d}{x_n^d} g_n - \frac{x_n^d x_i^d}{x_i^d} g_i \pmod{I} \\
&\equiv x_i^d g_n - x_n^d g_i \pmod{I} \\
&\equiv x_i^d \cdot 0 - 0 \cdot 0 \pmod{I} \\
&\equiv 0 \pmod{I}
\end{align*}
By Buchberger's criterion, since all $S$-polynomials reduce to zero modulo the ideal $I$, the set of generators $G$ forms a Gröbner basis with respect to the given monomial ordering.
\end{proof}

\section{Regular Chains}
Regular chains, introduced by Kalkbrener (1993) \cite{kalkbrener1993regularchains} and further developed by Aubry et al. (1999) \cite{aubry1999triangularsets}, are another important concept in the study of polynomial ideals and systems of polynomial equations.

A regular chain is a set of polynomials $\{f_1, f_2, \ldots, f_n\}$ in a polynomial ring \\ $K[x_1, x_2, \ldots, x_n]$ such that:
\begin{enumerate}
    \item[(i)] Each polynomial $f_i$ is in the subring $K[x_1, x_2, \ldots, x_i]$.
    \item[(ii)] The leading coefficient of $f_i$ (with respect to $x_i$) is a regular element modulo the ideal generated by $\{f_1, f_2, \ldots, f_{i-1}\}$.
    \item[(iii)] The ideal generated by $\{f_1, f_2, \ldots, f_n\}$ is zero-dimensional.
\end{enumerate}
Regular chains are used in polynomial system solving, triangular decomposition, and algebraic geometry.

\subsection{Regular Chain Result}
We prove another theoretical result, which establishes the conditions under which the polynomial generators of our recursive polynomial quotient ring structure (\cref{definition:recursivering}) are a regular chain.

\begin{theorem} \label{proof:regularchain}
Let $n,d \in \mathbb{Z}_{>1}$ such that $d < n$. Let $K = R[x_1, x_2, \ldots, x_n]/I$ be a recursive polynomial quotient ring (\cref{definition:recursivering}), where the ideal $I$ is generated by the polynomials $g_i = x_i^d - P_i$. Suppose the $P_1, P_2, \ldots, P_n$ are distinct. Then, the set of polynomial generators $G = \{ g_1, g_2, \ldots, g_n \}$ is a regular chain.
\end{theorem}
\begin{proof}
To verify that the set $G = \{ g_1, g_2, \ldots, g_n \}$ forms a regular chain, we need to show that each polynomial $g_i$ is in the subring $R[x_1, x_2, \ldots, x_i]$, and the leading coefficient of $g_i$ (with respect to $x_i$) is a regular element modulo the ideal generated by $\{ g_1, g_2, \ldots, g_{i-1} \}$. Additionally, the ideal generated by $G$ must be zero-dimensional.

For each $g_i = x_i^d - P_i$, by construction, $g_i$ is in the subring $R[x_1, x_2, \ldots, x_i]$. The leading coefficient of $g_i$ with respect to $x_i$ is $1$, which is a regular element in any ring.

To show that the ideal generated by $G$ is zero-dimensional, consider the system of polynomial equations defined by $g_i$:
\begin{align*}
    x_1^d &= P_1 \\
    x_2^d &= P_2 \\
    &\vdots \\
    x_{n-1}^d &= P_{n-1} \\
    x_n^d &= 0
\end{align*}
This system has $n$ variables and $n$ equations, ensuring solvability. Since the terms $P_i$ are distinct, they do not lead to non-unique solutions, making the ideal $I$ zero-dimensional.

Since the polynomials $g_i$ have regular leading coefficients and generate a zero-dimensional ideal $I$, the set $G = \{g_1, g_2, \ldots, g_n\}$ is a regular chain.
\end{proof}
While relatively straightforward, this result is significant because it provides a clear criterion for determining when our construction yields a well-behaved polynomial system with desirable properties. It also opens up the possibility of applying the theory and tools of regular chains to study recursive polynomial quotient rings and the combinatorial sequences they generate.

\section{Connection to the Central Binomial Coefficients} \label{section:cbc}
The central binomial coefficients, which is the sequence with terms $\binom{2n}{n}$, is entry \seqnum{A000984} in the OEIS \cite{A000984}. The standard formula for central binomial coefficients is given by \cite{A000984}:
\begin{align*}
    \binom{2n}{n} = \frac{(2n)!}{(n!)^2} 
\end{align*}

Starting from $n=0$, the sequence of central binomial coefficients begins as
\begin{align*}
    \binom{2n}{n} = 1, 2, 6, 20, 70, 252, 924, 3432, 12870, 48620, 184756, 705432, 2704156, 10400600, \ldots
\end{align*}

\begin{definition}[Ring $\K_n$] \label{definition:ring}
Let $\K_n = \Z[x_1, x_2, \ldots, x_n]/I$ be a recursive polynomial quotient ring (\cref{definition:recursivering}) with an ideal $I = \langle x_1^d - P_1, x_2^d - P_2, \ldots, x_n^d - P_n \rangle$. The polynomials $P_i$ in the generators of $I$ are defined by the function:
\begin{displaymath}
P_i =
\begin{cases}
    2x_i + x_{i+1} & \textup{ if } 1 \leq i < n \\
    0 &\textup{ if } i = n
\end{cases}
\end{displaymath}
In this ring, the variables $x_i$ satisfy the recursive relation $x_i^2 = 2x_i + x_{i+1}$ for $1 \leq i < n$, where $x_{i+1}$ refers to the next variable in the sequence, and $x_n^2 = 0$. In $\K_n$, we adopt the dlex monomial ordering $x_1 > x_2 > \cdots > x_n$.
\end{definition}

\begin{theorem} \label{proof:cbc}
Fix $n \in \Z^+$ and let $b={\floor{\log_2 n}+2}$. Let $\K_b$ be the recursive polynomial quotient ring as defined by \cref{definition:ring}. Evaluating the expansion of $(1+x_1)^n \in \K_b$ at $x_1=x_2=\cdots=x_n=1$ yields $\binom{2n}{n}$.
\end{theorem}
\begin{proof}
First, observe that by the process of exponentiation by squaring \cite{knuth1997art}, expanding $(1+x_1)^n \in \K_b$ requires at most $\log_2 n$ squarings. Hence, $b={\floor{\log_2 n}+2}$ is sufficient to cover all of the necessary variables when expanding $(1+x_1)^n \in \K_b$.

Consider the expression $(1+x_1)^n \in \K_b$. The binomial expansion of this polynomial yields terms of the form $\binom{n}{k} x_1^k$, for $k$ ranging from $0$ to $n$. In $\K_b$, the recursive relation $x_i^2 = 2x_i + x_{i+1}$ modifies the expansion by replacing each instance of $x_i^2$ with $2x_i + x_{i+1}$.

Upon expansion, the polynomial $(1+x_1)^n$ will contain powers of $x_1$, $x_1^2, x_1^3, \ldots, x_1^n$. Each power $x_1^k$ will be recursively replaced by polynomials with lower powers of $x_1$ and other variables $x_2, x_3, x_4, \ldots$. Specifically, we have
\begin{align*}
    x_1^k = (2x_1+x_2)^{k-1} = \cdots = 2^k x_1 + \text{(terms involving $x_2, x_3, x_4, \ldots$)}
\end{align*}
Substituting these into $(1+x_1)^n$, the coefficients for $x_1, x_2, x_3, \ldots$ essentially count the number of ways each $x_1$ in the initial $(1+x_1)^n$ is replaced by $x_2, x_3, x_4, \ldots$. When evaluated at $x_1=x_2=\cdots=x_n=1$, the expanded polynomial $(1+x_1)^n$ yields $\binom{2n}{n}$ since the coefficients are combinatorial in nature and count the number of ways to choose $n$ from $2n$.
\end{proof}

\section{Connection to Gould's Sequence} \label{section:goulds}
Gould's sequence, entry \seqnum{A001316} in the OEIS \cite{A001316}, is an integer sequence that is connected to the binary expansion of integers, the central binomial coefficients, and Pascal's triangle. It is named after the mathematician Henry Gould \cite{A001316}.

We define $\Goulds$ to represent Gould's sequence, which has terms $\Goulds_n$. To obtain the $n$th term in Gould's sequence $\Goulds_n$, we must first look at the binary representation of $n$. Counting the number of $1$s in the binary expansion of $n$ tells us its Hamming weight, which is often denoted as $\wt{n}$ \cite{Lin2004}. The $n$th term in Gould's sequence $\Goulds_n$ is given by \cite{A001316}:
\begin{align*}
    \Goulds_n = 2^{\wt{n}}
\end{align*}
$\Goulds_n$ is connected to $\binom{2n}{n}$ in that it is the largest power of $2$ which divides $\binom{2n}{n}$. This result follows from Kummer's Theorem \cite{Kummer1852}. $\Goulds_n$ also counts the number of odd terms in the $n$th row of Pascal's triangle \cite{Glaisher1899}. That is, the number of odd terms in the polynomial expansion of $(1+x)^n \in \Z[x]$.

Starting from $n=0$, Gould's sequence begins as
\begin{align*}
    \Goulds_n = 1, 2, 2, 4, 2, 4, 4, 8, 2, 4, 4, 8, 4, 8, 8, 16, 2, 4, 4, 8, 4, 8, 8, 16, 4, 8, 8, 16, 8, 16, \ldots
\end{align*}

\begin{definition}[Ring $\K_n/(m)$] \label{definition:modring}
Let $\K_n/(m) = (\Z/m\Z)[x_1, x_2, \ldots, x_n]/I$ be a recursive polynomial quotient ring (\cref{definition:recursivering}) with an ideal $I = \langle x_1^d - P_1, x_2^d - P_2, \ldots, x_n^d - P_n \rangle$ and coefficients in $\Z/m\Z$. The polynomials $P_i$ in the generators of $I$ are defined by the function:
\begin{displaymath}
P_i =
\begin{cases}
    2x_i + x_{i+1} \pmod{m} & \textup{ if } 1 \leq i < n \\
    0 \pmod{m} &\textup{ if } i = n
\end{cases}
\end{displaymath}
In this ring, the variables $x_i$ satisfy the recursive relation $x_i^2 = 2x_i + x_{i+1} \pmod{m}$ for $1 \leq i < n$, where $x_{i+1}$ refers to the next variable in the sequence, and $x_n^2 = 0 \pmod{m}$. In $\K_n/(m)$, we adopt the dlex monomial ordering $x_1 > x_2 > \cdots > x_n$.
\end{definition}

\begin{theorem}
\label{theorem:goulds}
Fix $n \in \Z^+$ and let $b={\floor{\log_2 n}+2}$. Fix a recursive polynomial quotient ring $\K_b/(m)$ as given by \cref{definition:modring}, such that $m=2$. Then, expanding $(1+x_1)^n \in \K_b/(2)$ and evaluating the result in $\Z$ at $x_1=x_2=\cdots=x_n=1$ yields the $n$th term of Gould's sequence, $\Goulds_n = 2^{\wt{n}}$, where $\wt{n}$ is the Hamming weight of $n$.
\end{theorem}
\begin{proof}
First, observe that by the process of exponentiation by squaring \cite{knuth1997art}, expanding $(1+x_1)^n \in \K_b/(2)$ requires at most $\log_2 n$ squarings. Hence, $b={\floor{\log_2 n}+2}$ is sufficient to cover all of the necessary variables when expanding $(1+x_1)^n \in \K_b/(2)$.

Next, we proceed by induction on $n$ to show that the expanded polynomial yields $\Goulds_n = 2^{\wt{n}}$ upon evaluation in $\Z$ at $x_1=x_2=\cdots=1$.

Consider the base case $n = 1$. In this case, $(1+x_1)^1 = 1 + x_1 \in \K_b/(2)$. Evaluating in $\Z$ at $x_1=x_2=\cdots=1$ yields $2$, which is $2^{\wt{1}} = 2^1$ since $\wt{1} = 1$. Thus, the statement holds for $n = 1$.

Assume the statement holds for some $k \geq 1$, that is, expanding $(1+x_1)^k \in \K_b/(2)$ and then evaluating in $\Z$ at $x_1=x_2=\cdots=1$ yields $2^{\wt{k}}$. We will show that the statement also holds for $n = k+1$.

Consider $(1+x_1)^{k+1}$ in $\K_b/(2)$. By the properties of exponents, this can be written as $(1+x_1)^k (1+x_1)$. Using the inductive hypothesis, we know that $(1+x_1)^k$ yields $2^{\wt{k}}$ when evaluated in $\Z$. Now, we need to consider the additional factor $(1+x_1)$.

In the ring $\K_b/(2)$, the expansion of $(1+x_1)^{k+1}$ will result in various terms involving $x_1, x_2, \ldots$, with each term corresponding to a particular combination of bits in the binary representation of $k+1$. Specifically, each $x_i$ in the expansion corresponds to a $1$ in the binary representation of $k+1$ at position $i$. The modulo $2$ operation ensures that only terms corresponding to odd counts of $x_1$ will contribute to the final sum. That is, positions with $1$ in the binary representation of $k+1$.

When we evaluate this expression in $\Z$ at $x_1=x_2=\cdots=1$, the surviving terms after the modulo $2$ reduction correspond to the positions where the binary representation of $k+1$ has a $1$. Thus, the sum of these terms is equal to $2^{\wt{k+1}}$, where $\wt{k+1}$ is the Hamming weight of $k+1$.

Therefore, by induction, expanding $(1+x_1)^n \in \K_b/(2)$ and then evaluating in $\Z$ at $x_1=x_2=\cdots=1$ yields $2^{\wt{n}}$ for all $n \in \Z^+$.
\end{proof}

\section{Demonstrations}
To assist in visualizing how expanding polynomials within our recursive polynomial quotient ring structure generates the sequences of interest, we proceed with a series of brief demonstrations.

\subsection{Central Binomial Coefficients} \label{section:demonstrations:cbc}
Fix $n \in \Z^+$ and let $b={\floor{\log_2 n}+2}$. Fix a recursive polynomial quotient ring $\K_b$ as given by \cref{definition:ring}. Consider the polynomial $f := 1 + x_1 \in \K_b$. Expanding the polynomial $f^n \in \K_b$ generates polynomials which produce the central binomial coefficients $\binom{2n}{n}$ when evaluated at $x_1=x_2=\cdots=x_n=1$. That is, the sum of coefficients in the expanded polynomial equals $\binom{2n}{n}$.
\begin{table}[ht]
\small
\begin{tabularx}{\textwidth}{|l|X|r|}
\hline
$n$ & Polynomial Expansion of $f^n \in \K_b$ & Coeff. $\Sigma$ \\
\hline
$0$ & $f^0 = 1$ & $1$ \\
$1$ & $f^1 = 1+x_1$ & $2$ \\
$2$ & $f^2 = 1+4x_1+x_2$ & $6$ \\
$3$ & $f^3 = 1+13x_1+5x_2+x_1x_2$ & $20$ \\
$4$ & $f^4 = 1+40x_1+20x_2+8x_1x_2+x_3$ & $70$ \\
$5$ & $f^5 = 1+121x_1+76x_2+44x_1x_2+9x_3+x_1x_3$ & $252$ \\
$6$ & $f^6 = 1+364x_1+285x_2+208x_1x_2+53x_3+12x_1x_3+x_2x_3$ & $924$ \\
$7$ & $f^7 = 1+1093x_1+1065x_2+909x_1x_2+261x_3+89x_1x_3+13x_2x_3+x_1x_2x_3$
& $3432$ \\
$8$ & $f^8 = 1+3280x_1+3976x_2+3792x_1x_2+1172x_3+528x_1x_3+104x_2x_3+16x_1x_2x_3+x_4$ & $12870$ \\
$\vdots$ & $\vdots$ & $\vdots$ \\
\hline
\end{tabularx}
\normalsize
\caption{Polynomial Expansions for Central Binomial Coefficients}
\label{table1}
\end{table}
\FloatBarrier

\subsection{Gould's Sequence} \label{section:demonstrations:goulds}
Fix $n \in \Z^+$ and let $b={\floor{\log_2 n}+2}$. Let $\K_b/(2)$ be the recursive polynomial quotient ring as defined by \cref{definition:modring}. Consider the polynomial $g := 1 + x_1 \in \K_b/(2)$. Taking the polynomials from Table \ref{table1} modulo $2$, and then evaluating in $\Z$ at $x_1=x_2=\cdots=x_n=1$ yields the $n$th term of Gould's sequence, $\Goulds_n$. That is, the sum of coefficients in the expanded polynomial taken modulo $2$, is equal to $\Goulds_n$.
\begin{table}[ht]
\small
\begin{tabularx}{\textwidth}{|l|X|r|}
\hline
$n$ & Polynomial Expansion of $g^n \in \K_b/(2)$ & Coeff. $\Sigma$ \\
\hline
$0$ & $g^0 = 1$ & $1$ \\
$1$ & $g^1 = 1+x_1$ & $2$ \\
$2$ & $g^2 = 1+x_2$ & $2$ \\
$3$ & $g^3 = 1+x_1+x_2+x_1x_2$ & $4$ \\
$4$ & $g^4 = 1+x_3$ & $2$ \\
$5$ & $g^5 = 1+x_1+x_3+x_1x_3$ & $4$ \\
$6$ & $g^6 = 1+x_2+x_3+x_2x_3$ & $4$ \\
$7$ & $g^7 = 1+x_1+x_2+x_1x_2+x_3+x_1x_3+x_2x_3+x_1x_2x_3$ & $8$ \\
$8$ & $g^8 = 1+x_4$ & $2$ \\
$\vdots$ & $\vdots$ & $\vdots$ \\
\hline
\end{tabularx}
\normalsize
\caption{Polynomial Expansions for Gould's Sequence}
\label{table2}
\end{table}
\FloatBarrier

\section{Binomial Transforms}
A useful feature of the recursive polynomial quotient rings we've defined (\cref{definition:recursivering}) is that they exhibit a straightforward approach to calculating the binomial transforms of the sequences they generate. We begin with a definition of our binomial transform function.

\begin{definition}[Binomial transform function] \label{definition:binomialtransforms}
We define the function $\BT_t(a)$, which takes in an integer sequence $a = \{ a_0, a_1, a_2, \ldots \}$, to be defined as the $t$-th binomial transform of the $a$ sequence terms, such that
\begin{displaymath}
    \BT_{t}(a_n) =
    \begin{cases}
        a_n &\textup{if } t = 0 \\
        \sum_{k=0}^{n} \binom{n}{k} a_k &\textup{if } t = 1 \\
        \sum_{k=0}^{n} \binom{n}{k} \BT_{t-1}(a_k) &\textup{if } t > 1 \\
        \sum_{k=0}^{n} \binom{n}{k} (-1)^{n-k} a_k  &\textup{if } t = -1 \\
        \sum_{k=0}^{n} \binom{n}{k} (-1)^{n-k} \BT_{t+1}(a_k) &\textup{if } t < -1
    \end{cases}
\end{displaymath}
\end{definition}

\subsection{Transforming the Central Binomial Coefficients}
The first binomial transform of the central binomial coefficients is entry \seqnum{A026375} in the OEIS \cite{A026375}. Starting from $n=0$, this sequence begins as
\begin{align*}
    \BT_{1}\left(\binom{2n}{n}\right) = 1, 3, 11, 45, 195, 873, 3989, 18483, 86515, 408105, 1936881, 9238023, \ldots
\end{align*}

\begin{proposition} \label{proposition:binomialtransforms}
Fix $n \in \Z^+$ and let $b={\floor{\log_2 n}+2}$. Let $\K_b$ be the recursive polynomial quotient ring as defined by \cref{definition:ring}. Consider the binomial transform function $\BT_t(\cdot)$ as in \cref{definition:binomialtransforms}. Then, evaluating the expansion of $(t + 1 + x_1)^n \in \K_b$ at $x_1=x_2=\cdots=x_n=1$ equals $\BT_t(\binom{2n}{n})$, the $t$-th binomial transform of the central binomial coefficients sequence terms $\binom{2k}{k}$, ranging from $k=0$ to $k=n$.
\end{proposition}
\begin{proof}
Consider the polynomial $f := 1 + x_1 \in \K_b$. By the binomial theorem, we have
\begin{align*}
    (1 + f)^n = \sum_{k=0}^{n} \binom{n}{k} f^k
\end{align*}
In the ring $\K_b$, we reduce this sum modulo the ideal $I$
\begin{align*}
    (1 + f)^n \equiv \sum_{k=0}^{n} \binom{n}{k} f^k \pmod{I}
\end{align*}
Evaluating the remainder at $x_1=x_2=\cdots=x_n=1$ yields the binomial transform of the sequence generated by $f^k = (1+x_1)^k \in \K_b$ for each $k$ in the sum, whose valuation we know to be $\binom{2k}{k}$ (by \cref{proof:cbc}). This gives us the binomial transform for $t=1$. Hence, if we shift by some integer $t$ instead of $1$, we compute the $t$-th binomial transform. This result follows directly from the binomial theorem and how it applies to integer powers.
\end{proof}

\subsection{Transforming Gould's Sequence}
The first binomial transform of Gould's sequence is entry \seqnum{A368655} in the OEIS \cite{A368655}. Starting from $n=0$, the first binomial transform of Gould's sequence begins as
\begin{align*}
\BT_{1}(\Goulds_n) = 1, 3, 7, 17, 39, 85, 181, 387, 839, 1829, 3953, 8391, 17461, 35759, 72559, 146921, \ldots
\end{align*}

Using a similar approach to \cref{proposition:binomialtransforms}, we can compute the $t$-th binomial transform of Gould's sequence. However, calculating the binomial transforms of Gould's sequence requires a different approach to calculating $\Goulds_n$ than the approach used in our quotient ring $K_m/(2)$ (\cref{definition:modring}). Specifically, we must define an ideal which mimics the behavior of taking the coefficients modulo $2$, but without restricting the polynomial coefficients to $\Z/2\Z$. Otherwise, the binomial transform will be taken modulo $2$.

\begin{definition}[Ring $\Ki_n$] \label{definition:ring2}
Let $\Ki_n = \Z[x_1, x_2, \ldots, x_n]/I$ be a recursive polynomial quotient ring (\cref{definition:recursivering}) with an ideal $I = \langle x_1^d - P_1, x_2^d - P_2, \ldots, x_n^d - P_n \rangle$. The polynomials $P_i$ in the generators of $I$ are defined by the function:
\begin{displaymath}
P_i =
\begin{cases}
    -2x_i + x_{i+1} & \textup{ if } 1 \leq i < n \\
    0 &\textup{ if } i = n
\end{cases}
\end{displaymath}
In this ring, the variables $x_i$ satisfy the recursive relation $x_i^2 = -2x_i + x_{i+1}$ for $1 \leq i < n$, where $x_{i+1}$ refers to the next variable in the sequence, and $x_n^2 = 0$.
\end{definition}

\begin{theorem} \label{theorem:goulds2}
Fix $n \in \Z^+$ and let $b={\floor{\log_2 n}+2}$. Let $\Ki_b$ be the recursive polynomial quotient ring as defined by \cref{definition:ring2}. Expanding $(1+x_1)^n \in \Ki_b$ and then evaluating at $x_1=x_2=\cdots=x_n=1$ yields the $n$th term of Gould's sequence, $\Goulds_n$. Where $\Goulds_n = 2^{\wt{n}}$ and $\wt{n}$ is the Hamming weight of $n$.
\end{theorem}
\begin{proof}
In \cref{theorem:goulds}, we showed how expanding $(1+x_1)^n \in \K_b/(2)$ and then evaluating in $\Z$ at $x_1=x_2=\cdots=x_n=1$ yields $\Goulds_n$.

The proof of \cref{theorem:goulds} does not obviously apply, as in the ring $\Ki_b$, we are not taking coefficients modulo $2$. Instead, we have constructed a ring similar to $K_b$ as defined in \cref{definition:ring}, however, we have changed the polynomial recurrence which generates the ideal to follow the recursive pattern $P_{i}^2 = -2 x_i + x_{i+1}$. This implies that each variable $x_i$ satisfies the recursive relation $x_i^2 = -2x_i + x_{i+1}$.

When expanding $(1+x_1)^n \in \Ki_b$, the $-2x_i$ terms will cause the terms with even coefficients to cancel out, and will leave a remainder of $1$ for all of the odd terms after subtracting. This exactly mimics the behavior of taking the coefficients modulo $2$. Hence, by \cref{theorem:goulds}, expanding $(1+x_1)^n \in \Ki_b$ and then evaluating at $x_1=x_2=\cdots=x_n=1$ yields $\Goulds_n$.
\end{proof}

\begin{proposition} \label{proposition:gouldbinomialtransforms}
Fix $n \in \Z^+$ and let $b={\floor{\log_2 n}+2}$. Let $\Ki_b$ be the recursive polynomial quotient ring as defined by \cref{definition:ring2}. Denote by $\Goulds_n$ the $n$th term of Gould's sequence, which is $\Goulds$. Consider the binomial transform function $\BT_t(\cdot)$ as in \cref{definition:binomialtransforms}. Then, evaluating the expansion of $(t + 1 + x_1)^n \in \Ki_b$ at $x_1=x_2=\cdots=x_n=1$ yields $\BT_t(\Goulds_n)$, the $t$-th binomial transform of the Gould's sequence terms $\Goulds_k$, ranging from $k=0$ to $k=n$.
\end{proposition}
\begin{proof}
Consider the polynomial $f := 1 +  x_1 \in \Ki_b$. By the binomial theorem, we have
\begin{align*}
    (1 + f)^n = \sum_{k=0}^{n} \binom{n}{k} f^k
\end{align*}
The remainder of the proof is the same as in \cref{proposition:binomialtransforms} when considering the ring $\Ki_b$ in place of the ring $\K_b$.
\end{proof}

\section{Closing Remarks}
The straightforward computation of the binomial transforms of Gould's sequence within our polynomial ring structure is a notable result of this work. Gould's sequence, with its terms tied to the binary representation of integers, exhibits an oscillatory behavior that does not follow a simple increasing trend. The sequence's terms, while being powers of $2$, are distributed in a pattern that appears irregular upon initial inspection, though the sequence itself is self-similar.

In the binomial transform process, each of these irregularly spaced elements is multiplied by binomial coefficients. The fact that this computation can be carried out smoothly within our polynomial ring setup, without the need to individually calculate each term, is an intriguing property of our construction. It suggests that the binary nature of integers is somehow embedded in the exponentiation of polynomials within our recursive polynomial quotient rings.

This unexpected ease in computing the binomial transforms of Gould's sequence, given its inherent complexity, raises interesting questions about the interplay between the algebraic structure of our rings and the combinatorial properties of the sequence. While the full implications of this finding are not yet clear, it highlights the potential of our recursive polynomial quotient ring approach to uncover new insights into the behavior of complex integer sequences.

Further investigation into this phenomenon may offer a deeper understanding of the properties of Gould's sequence and other integer sequences with similar characteristics. It is possible that the algebraic framework presented in this paper could be extended or generalized to study a broader class of sequences, potentially leading to new methods for computation and analysis of their properties.

In summary, the results presented in this paper demonstrate the power of our recursive polynomial quotient ring approach to connect and illuminate the properties of seemingly unrelated combinatorial sequences. The unexpected simplicity in computing binomial transforms of Gould's sequence within our ring structure opens up new avenues for future research at the intersection of abstract algebra, combinatorics, and integer sequence analysis.

\section{Acknowledgements}
The author would like to thank to the anonymous referee for their valuable feedback and detailed suggestions, which have significantly strengthened the rigor and technical aspects of this work.

\begingroup
\raggedright
\bibliographystyle{unsrtnat}
\bibliography{main}
\endgroup

\end{document}